\documentclass[]{llncs}
\usepackage{amsfonts}
\usepackage{amsmath}
\usepackage{amssymb}
\usepackage{graphics}
\usepackage[dvips]{epsfig}
\usepackage{times}
\usepackage{algorithm}
\usepackage{algorithmic}
\usepackage{fancyhdr}
\usepackage{ifthen}
\newcommand\settitle[2][]{%
 \title{#2}
 \ifthenelse{\equal{#1}{}}%
  {\fancyhead[RO]{\nouppercase #2 \qquad \thepage}}%
  {\fancyhead[RO]{\nouppercase #1 \qquad \thepage}}%
}

\fancypagestyle{plain}{%
\fancyhf{}

}

   \newtheorem{ex}{Example}
   \newtheorem{Bemerkung}{Remark}

   \newcommand{\D}{\displaystyle}

   \def\C{{\mathbb C}}
   \def\R{{\mathbb R}}

   \def\Pf{{\it Proof.$\;\;$}}

   \def\qed{\hfill$\diamond$}

   \def\cP{{\mathcal P}}

\def\C{{\mathbb C}}

\def\R{{\mathbb R}}

\def\Pf{{\it Proof.$\;$}}

\def\qed{\hfill$\blacksquare$}

\def\lin{\mbox{\rm lin}}
\def\Pr{\mbox{\rm Pr}}

\def\({\langle}
\def\){\rangle}
\def\mb{\boldsymbol}

\def\cF{{\mathcal F}}

\def\cH{{\mathcal H}}

\def\cP{{\mathcal P}}

\def\cU{{\mathcal U}}

\def\1{\mb1}

\def\v0{{\bf 0}}

\def\ov{\overline}

\begin{document}

\settitle[Markovian Statistics]{Markovian Statistics on Evolving Systems}
\author{Ulrich Faigle ${}^1$, Gerhard Gierz ${}^2$}
\institute{Mathematisches Institut\\
           Universt\"at zu K\"oln\\
           Weyertal 80, 50931 K\"oln, Germany\\
\email{faigle@zpr.uni-koeln.de}\\[1ex]
\and
Department of Mathematics\\
University of California at Riverside\\
Riverside, CA 92521, USA\\
\email{gierz@math.ucr.edu}}

\date{}
\maketitle

\thispagestyle{plain}
\begin{abstract}
A novel framework for the analysis of observation statistics on time discrete linear evolutions in Banach space is presented. The model differs from traditional models for stochastic processes and, in particular, clearly distinguishes between the deterministic evolution of a system and the stochastic nature of observations on the evolving system. General Markov chains are defined in this context and it is shown how typical traditional models of classical or quantum random walks and Markov processes fit into the framework and how a theory of quantum statistics ({\it sensu} Barndorff-Nielsen, Gill and Jupp) may be developed from it. The framework permits a general theory of joint observability of two or more observation variables which may be viewed as an extension of the Heisenberg uncertainty principle and, in particular, offers a novel mathematical perspective on the violation of Bell's inequalities in quantum models. Main results include a general sampling theorem relative to Riesz evolution operators in the spirit of von Neumann's mean ergodic theorem for normal operators in Hilbert space.
\end{abstract}

\begin{keywords} {Banach space, Bell inequality, ergodic theorem, evolution, Heisenberg uncertainty, Hilbert space, Markov chain, observable, quantum statistics, random walk, sampling, Riesz operator, stochastic process}
\end{keywords}

\section{Introduction}\label{sec:Introduction}
Consider a system $\mathfrak S$ that is observed at discrete times $t=1,2,3,\ldots$ relative to a pre-specified event $S$ that may or may not occur at time $t$. A statistical analysis is interested in the relative frequency of the event $S$. In particular, one would like to know whether the sample frequency of its occurrence converges to a definite limiting value. For a mathematical formulation of this problem, let us define the indicator functions $I_S^{(t)}$ as ($0,1$)-variables with the event notation $\{I^{(t)}_S=1\}$ meaning that $S$ is observed at time $t$. The sampled frequency up to  time $t$ would then be
$$
      \ov{I}^{(t)}_S = \frac1t\sum_{m=1} I^{(m)}_S.
$$
We call the study of the limiting behavior of the sample frequency of an event $S$ as the \emph{Markovian} problem.

\medskip
Mathematical approaches to this problem typically think of $\mathfrak S$ as a random source that emits symbols from some (finite or infinite) alphabet $A$ and thus gives rise to an  $A$-valued stochastic process $(X_t)$. $S$ is assumed to represent a particular event relative to this process that may or may not materialize at time $t$. So the observation variables $I^{(t)}_S$ reflect an associated stochastic process in their own right and the expected average number of observations of $S$ at time $t$ is
$$
    E(\ov{I}^{(t)}_S) = \frac1t\sum_{m=1}^t E(I^{(m)}_S) = \frac1t\sum_{m=1}^t I^{(m)}\Pr\{I^{(m)}_S=1\}.
$$
No matter what the nature of the underlying stochastic process $(X_t)$ is, the mathematical analysis of the statistical problem relative to the event $S$, \emph{i.e.}, the Markovian problem, will actually be on the associated binary process $(I^{(t)}_S)$.

\medskip
The present investigation is concerned with the Markovian problem of the relative occurrence of some event $S$ rather than with the mathematical analysis of general stochastic processes {\it per se}. Therefore, there is no loss in generality when we restrict ourselves to a model where $\mathfrak S$ is viewed as some source that emits symbols from some alphabet $A$ of \emph{finite} cardinality $|A|<\infty$. In fact, the assumption of $A$ as a binary alphabet (\emph{i.e.}, $|A|=2$) would theoretically suffice. However, it is convenient to also consider more general alphabets occasionally.

\medskip
A stochastic process $(X_{t})$ is usually understood
to be a sequence of stochastic variables $X_{t}$ that are defined on some
probability space and one is interested in the expected limiting
behavior of $(X_{t})$. However, if one takes statistics on $(X_{t})$,
\emph{i.e.}, averages the observation of special events over time, there is not always
a clear asymptotic behavior. Indeed, there are examples of even completely deterministic
processes $(X_{t})$ where observation statistics do not converge. It is the limiting behavior
of observation averages we are concerned with here.

\medskip
Our approach is motivated by a Markovian interpretation of $(X_{t})$:
A random source ${\mathfrak{S}}$ produces symbols $a$ of an alphabet $A$ over
time $t$ with $\{X_{t}=a\}$ de\-noting these events. ${\mathfrak{S}}$ is thought
to change over discrete time $t$, with the probability
$\mbox {\rm Pr}\{X_{t}=a\}$ depending on the current state of ${\mathfrak{S}}%
$. The classical example goes back to Markov's~\cite{Markoff} model of a
system that admits a set $N$ of ground states and is subject to a random walk
on $N$ with transition probabilities $p_{ij}$. The system states are then the
probability distributions $p^{(t)}$ of the positions of the random walk at times $t$. The
Markov model has been very successful in application modeling. Statistical
mechanics in physics, for example, describes the behavior of ideal gases in this
way. But also the behavior of economic and social systems is often viewed as
following Markovian principles. Internet search engines successfully organize and
rank their search according to Markovian statistics.

\medskip The situation seems to be more complicated with quantum systems that
do not admit a classical analysis. For example, the result of a quantum measurement
is not a deterministic function of the state of the system and the measuring instrument
applied but rather an expected value relative to some (state dependent) probability distribution
on the possible measurement outcomes. Moreover, Heisenberg's uncertainty principle says that
observations may not be simultaneously feasible unless they conform to a
special condition. Experimental evidence with spin correlations
(Aspect \textit{et al.}~\cite{Aspect}) furthermore exhibits a definite violation of classical
statistical principles as expressed in Bell's~\cite{Bell64,Bell66}
inequalities. While the Schr\"odinger picture of quantum states being
described by wave functions yields a special theory of quantum probabilities
with applications also to quantum computing, active current research effort is
devoted to the quantum analogs of classical Markov random walks. In this spirit,
Barndorff-Nielsen, Gill and Jupp~\cite{B-NGJ} have put forward a theory of
quantum statistical inference.

\medskip
The present investigation proposes a quite general model for Markovian statistical analysis. Rather than following the standard approach to stochastic processes, our model is linear and motivated by the linear algebraic analysis of classical Markov chains of Gilbert~\cite{Gilbert}, Dharmakhari~\cite{Dharma65} and Heller~\cite{Heller}, which has led to the identification of more general Markov type processes (\emph{e.g.}, Jaeger~\cite{Jaeger}). Addressing the issue of the ''dimension'' of a stochastic evolution, the asymptotic behavior of even more general stochastic processes could be clarified (Faigle and Sch\"onhuth~\cite{FS}). Generalizing these previous models, our setting is in Banach space and focusses on the evolution of linear operators, which allows us to deal also with the statistics of discrete quantum type evolutions appropriately.

\medskip
There are several advantages and novel aspects in our approach. Not only does our model include typical Markovian models proposed so far (see the examples in Section~\ref{sec:examples}), but its generality allows us to develop a meaningful notion of \emph{jointly observable} statistical measuring instruments on an evolving system. Sets of classical stochastic variables are always jointly observable (for the simple reason that they are mathematically based on the same underlying probability space). The Heisenberg uncertainty principle, on the other hand, makes it clear that this property is no longer guaranteed for statistical observations on quantum systems.

\medskip
While the Heisenberg principle is formulated for pairs of self-adjoint operators, our model allows us to deal with three (or more) operators as well. We show that the Heisenberg principle corresponds to a very special case in our setting (Section~\ref{sec:joint-observations}). In fact, a careful mathematical analysis of the joint observability of $3$ measurement operators may offer a straightforward key to the understanding of Bell's inequalities (Section~\ref{sec:Bell}).

\medskip
These advantages are the result of a clear separation of the aspect of the (deterministic) evolution of a system from the aspect of statistical observations on the evolving system in the mathematical model.

\medskip
Our presentation is organized as follows. Section~\ref{sec:evolutions} introduces evolution operators on Banach spaces and discusses their ergodicity. Then sampling functions are studied and their convergence behavior is characterized in the Sampling Theorem (Theorem~\ref{t.sampling-theorem}) relative to finitary  evolutions, which include all evolutions based on Riesz operators, for example. Observables and generalized Markov chains are defined in Section~\ref{sec:M-chains}. These notions are illustrated by the examples in Section~\ref{sec:examples} with particular emphasis on random walks and quantum statistics. The proofs of the main results are deferred into the Appendix.

\section{Evolutions of systems}\label{sec:evolutions} Let ${\mathfrak{S}}$ be some \emph{system} that is in a certain \emph{state} $S_t$ at any time $t$. Observing $\mathfrak S$ at discrete times $t=0,1,2,\ldots$, we refer to the sequence $\epsilon =(S_t)_{t\geq 0}$ as an \emph{evolution} of $\mathfrak S$. For a mathematical analysis, the evolution needs to be represented in some (mathematical) \emph{universe} $\cU$. In the present investigation, we will always assume $\cU$ to be be a vector space over the complex field $\C$. A  \emph{representation} of the evolution $\epsilon$ in $\cU$ is then a map $t\mapsto s^{(t)}\in \cU$ such that there is a linear operator $\psi$ on $\cU$ with the property
$$
       s^{(t+1)} =  \psi s^{(t)} \quad(t=0,1,2,\ldots).
$$
We think of the vector $s^{(t)}\in \cU$ as the representation of the state $S_t$ of $\mathfrak S$ at time $t$ and call $\psi$ an \emph{evolution operator}. Clearly, any evolution $\epsilon$ of $\mathfrak S$ admits such a representation. For example, a \emph{stationary} representation, where $s^{(t)} =s^{(0)}$ for all $t$ and the evolution operators are exactly those (linear) operators on $\cU$ that fix $s^{(0)}$.

\medskip
In a practical system analysis, it is the first task of the modeler consists in the  determination of an appropriate representation of the evolution of the system $\mathfrak S$ under consideration. Here, however, we will assume that the evolution is already represented in some universe $\cU$ so that the evolutions are vector sequences $\Psi$ of the form
$$
    \Psi=(\psi, s) = (\psi^t s\mid t=0,1,2,\ldots)
$$
where $\psi$ is an operator on $\cU$. We furthermore assume that $\cU$ is endowed with some norm $\|\cdot\|$  and is complete with respect to this norm (otherwise we replace $\cU$ by its completion $\ov{\cU}$).

\begin{Bemerkung} By standard complexification arguments in functional analysis (\emph{e.g.}, \cite{Conway90}), the results we obtain in this section apply  to universes over the real field $\R$ as well. We choose $\C$ for mathematical con\-venience, without loss of generality.
\end{Bemerkung}

\medskip
The \emph{evolution space} of the evolution $\Psi=(\psi,s)$ in $\cU$ is the linear subspace $\cU_\Psi$ generated by $\Psi$, \emph{i.e.},
$$
{{\mathcal{U}}}_{\Psi}=\mbox {\rm lin}\{\psi^{t}s\mid t=0,1,\ldots \}.
$$
The parameter $\dim \Psi= \dim\cU_\Psi$ is the \emph{dimension} of the evolution $\Psi$. We will refer to the vectors  $s^{(t)}=\psi^{t}s$ as the \emph{states} of $\Psi$.

\medskip
Notice that $\cU_\Psi$ is $\psi$-invariant (\emph{i.e.}, $\psi(U_\Psi)\subseteq U_\Psi$).
So the restriction of $\psi$ to $U_\Psi$ is an operator
on the normed space $U_\Psi$. Let $\ov{\cU}_\Psi$ be the closure of $U_\Psi$ in $\cU$ and recall from general operator theory\footnote{\emph{e.g.}, \cite{Conway90,Dowson78}} that $\psi$ extends to a unique norm bounded (and hence continuous) operator $\ov{\psi}: \ov{\cU}_\Psi \to \ov{\cU}_\Psi$ with the same (finite) norm, provided $\psi$ is norm bounded on $\cU_\Psi$. The norm of $\psi$ on $\cU_\Psi$ is
\[
\Vert\psi\Vert_{s}=\operatorname{inf}\{c\in{{\mathbb{R}}}\mid\Vert\psi
u\Vert\leq c\,\Vert u\Vert\;\mbox {for all $u\in {{\mathcal{U}}}_{\Psi }$}\}.
\]
In the case of a finite-dimensional evolution (\emph{i.e.}, $\operatorname{dim} \Psi<\infty$), for
example, $\Vert\psi\Vert_{s}$ is necessarily finite. The \emph{norm} of the
evolution $\Psi=(\psi,s)$ is defined as
\[
\Vert\Psi\Vert =\operatorname{inf}\{c\in{{\mathbb{R}}}\mid\Vert\psi
^{t}s\Vert\leq c\,\Vert s\Vert\;,\forall t\geq0\}
\]
and $\Psi$ said to be \emph{stable} if $\Vert\Psi\Vert <\infty$. Thus
$\Psi$ is stable if $\Vert\psi\Vert_{s}\leq1$, for example.

\begin{lemma}\label{l.stability} If $\Psi=(\psi,s)$ is stable, then $\Psi'=(\psi,s')$ is stable for every $s'\in \cU_\Psi$. Hence the restriction of $\psi$ to $\cU_\Psi$ does not admit any eigenvalue $\lambda$ with $|\lambda| > 1$.
\end{lemma}

\Pf Consider any $s' =\sum_{j=1}^k a_j \psi^{t_j}s\in \cU_\Psi$. Then the triangle inequality yields
$$
\|\psi^t s'\| \leq \big(\|\Psi\|\sum_{j=1}^k |a_j|\big)\|s\| \quad\mbox{for all $t\geq 0$.}
$$
\qed

\medskip
The evolution $\Psi=(\psi,s)$ is \emph{ergodic} if
its states $s^{(t)}=\psi^{t}s$ converge in the norm. $\Psi$ is \emph{mean ergodic} if the state averages
\[
\overline{s}^{(t)}=\frac{1}{t}\sum_{m=1}^{t}s^{(m)}%
\]
converge to some \emph{limit state} $\overline{s}^{\infty} \in \ov{U}_\Psi$. Clearly, if $\ov{s}^{(\infty)}$ exists (if and $\|\psi\|_s<\infty$ holds), $\Psi$ is \emph{stationary} in the sense
$$
    \ov{\psi}\ov{s}^{(\infty)} = \ov{s}^{(\infty)}.
$$
Moreover, an ergodic evolution is also mean ergodic, while the converse conclusion is generally false.

\subsection{Equivalent evolutions}\label{sec:equivalent-evolutions}
Let us call two evolutions $\Phi=(\varphi,v)$ and  $\Psi=(\psi,w)$ in $\cU$ \emph{equivalent} if
$$
   \lim_{t\to\infty} \|\psi^t w - \varphi^t v\| = 0.
$$
By Cauchy's Theorem, equivalence implies
$$
\lim_{t\to\infty} \frac1t\sum_{m=1}^t\|\psi^m w - \varphi^m v\| = 0.
$$
So, assuming equivalence, $\Phi$ is mean ergodic exactly when $\Psi$ is mean ergodic and in either case, one has
\begin{equation}\label{eq.mean-ergodic-equivalence}
\lim_{t\to\infty} \frac1t\sum_{m=1}^t \psi^m w = \lim_{t\to\infty} \frac1t\sum_{m=1}^t\varphi^m v.
\end{equation}

\subsubsection{Finitary evolutions}\label{sec:finitary-evolutions}
We say that the evolution $\Psi=(\psi,s)$ is \emph{finitary} if $\Psi$ is equivalent to a finite-dimensional evolution $\Phi=(\varphi,v)$.  A characterization of the  mean ergodicity of a finite-dimensional evolution $\Phi$ follows from the analysis of Faigle and Sch\"onhuth~\cite{FS} and says in essence:  $\Phi$ is mean ergodic precisely when $\Phi$ is equivalent to an evolution $\Pi=(\pi,v)$, where $\pi$ is a projection operator on $\ov{\cU}_\Phi= \cU_\Phi$.

\begin{proposition}[\cite{FS}]\label{p.FS} Let $\Phi=(\varphi,v)$ be a finite-dimensional evolution. Then $\Phi$ is mean ergodic if and only if $\Phi$ is stable. Moreover, if $\Phi$ is stable, one has
$$
\lim_{t\to\infty} \frac1t\sum_{m=1}^t\varphi^m v =\left\{\begin{array}{cl} 0 &\mbox{if $\lambda = 1$ is not an eigenvalue of $\ov{\varphi}$}\\
P_1 v &\mbox{if $\lambda=1$ is an eigenvalue of $\ov{\varphi}$,}\end{array} \right.
$$
where  $P_1$ is a projection operator onto the eigenspace $E_1=\{x\in \cU_\Phi\mid x= \varphi x \}$.
\end{proposition}

\Pf By Lemma~\ref{l.stability}, $\Phi$ is stable if and only if $(\phi,v')$ is stable for all $v'\in \cU_\Phi$. So Proposition~\ref{p.FS} is a direct consequence of Theorem~2 and its proof in \cite{FS}.

\qed

\subsubsection{Riesz evolutions}\label{sec:Riesz-evolutions} Recall that the \emph{spectrum} $\sigma(T)$ of a (linear) operator $T:V\to V$ on a complex normed vector space $V$ consists of those $\lambda\in \C$ such that the operator $L_\lambda=T-\lambda$ (with values $L_\lambda v = Tv-\lambda v$) is not invertible. $T$ is called a \emph{Riesz operator} if $T$ is bounded and
\begin{enumerate}
\item[(a)] each $\lambda\in \sigma(T)\setminus\{0\}$ is an eigenvalue of $T$ with finite algebraic multiplicity;
\item[(b)] $0$ is the only possible accumulation point of $\sigma(T)$.
\end{enumerate}

Riesz operators form a quite wide class of operators that includes the so-called \emph{compact} operators.  In particular, every operator $T$ with finite-dimensional range is Riesz.
Further examples are the \emph{Hilbert-Schmidt} operators on a Hilbert space $\cH$, namely the bounded operators $T:\cH\to\cH$ such that
$$
    \sum_{i\in I} \|T e_i\|^2 < \infty
$$
holds for some orthonormal basis $\{e_i\mid i\in I\}$ of $\cH$. Note that every operator on a finite-dimensional vector space $V$ is trivially Hilbert-Schmidt relative to any inner product.

\medskip
A \emph{Riesz evolution} in our universe $\cU$ is now an evolution $\Psi=(\psi,s)$ such that $\psi$ extends to a Riesz operator $\ov{\psi}: \ov{\cU}_\Psi \to\ov{\cU}_\Psi$. In particular, every evolution under a Riesz evolution operator on $\cU$ is Riesz.

\medskip
The characterization of mean ergodic finite-dimensional evolutions (Proposition~\ref{p.FS}) extends to general Riesz evolutions.

\begin{theorem}\label{t.Riesz} Let $\Psi=(\psi,s)$ be any Riesz evolution in $\cU$. Then $\Psi$ is finitary. Moreover, $\Psi$ is mean ergodic if and only if $\Psi$ is stable.
In particular, if $\Psi$ is stable, one has
$$
\lim_{t\to\infty} \frac1t\sum_{m=1}^t\psi^m v =\left\{\begin{array}{cl} 0 &\mbox{if $\lambda = 1$ is not an eigenvalue of $\ov{\psi}$}\\
P_1 v &\mbox{if $\lambda=1$ is an eigenvalue of $\ov{\psi}$,}\end{array} \right.
$$
where  $P_1$ is a projection operator onto the eigenspace $E_1=\{x\in \hat{\cU}_\Phi\mid x= \ov{\psi} x \}$.
\end{theorem}

The essential part of the proof of Theorem~\ref{t.Riesz} consists in showing that Riesz evolutions are finitary. We discuss the details in the Appendix (\emph{cf.} Proposition~\ref{p.1} there).

\subsubsection{Normal evolutions in Hilbert space}\label{sec:normal-operators}
Let ${{\mathcal{H}}}$ be a Hilbert space and
recall that an operator $T$ on ${{\mathcal{H}}}$ is \emph{normal} if $T$
commutes with its adjoint $T^{*}$ (\emph{i.e.}, $TT^{*}=T^{*}T$).

\begin{theorem}\label{t.vNeumann-normal-mean} Let $\psi$ be a bounded normal operator on
${{\mathcal{H}}}$ and $\Psi=(\psi,s)$ an evolution. Then the following statements are equivalent:
\begin{enumerate}
\item[(i)] $\Psi$ is stable.
\item[(ii)] $\Psi$ is mean ergodic.
\end{enumerate}
In this case, the averages $\overline{s}^{(t)}$ converge to the ortho\-gonal projection
of $s$ onto the eigenspace $E_{1}=\{x\in{{\mathcal{H}}}\mid\psi
x=x\}$.

\end{theorem}

\medskip
We prove Theorem~\ref{t.vNeumann-normal-mean} in the Appendix. The implication ''(i) $\Rightarrow$ (ii)'' is well-known and usually stated as \emph{von Neumann's mean ergodic theorem}:

\begin{corollary}[von Neumann]\label{c.vNeumann-normal-mean} If $\psi$ is a normal operator on
${{\mathcal{H}}}$ of norm $\|\psi\|\leq 1$, then every evolution $\Psi=(\psi,s)$ is mean ergodic.
\qed
\end{corollary}

\medskip An important special case is the evolution of a wave function
$v\in{{\mathcal{H}}}$ of a quantum system in discrete time. According to
Schr\"odinger's differential equation, there is a unitary operator $U$
(\emph{i.e.} $UU^{*}=I=U^{*}U$) so that the discrete evolution of $v$ is given as
\[
v^{(t)}=U^{t}v\quad(t=0,1,\ldots).
\]
Clearly, the operator $v\mapsto Uv$ is normal and bounded. Moreover, $(U,v)$ is stable for any $v\in \cH$. So \emph{Schr\"odinger evolutions} are mean ergodic.

\subsection{Sampling}\label{sec:sampling}
By a \emph{sampling function} relative to the universe $\cU$ we understand a continuous linear map $f:{{\mathcal{U}}%
}\to \cF$, where $\cF$ is a normed vector space of \emph{samples}.  With respect to an evolution $\psi=(\psi,s)$, the $f_{t}=f(\psi^t s)$ are the
\emph{sampling values} with the corresponding \emph{sampling averages}
\[
\overline{f}_{t}=\frac{1}{t}\sum_{m=1}^{t}f_{m}=\frac{1}{t}\sum_{m=1}^{t}f(\psi^ms)  \quad(t=1,2,\ldots).
\]
 In applications, a sampling function $f$ will typically be a functional into the scalar field of $\cU$. But more general sample spaces $\cF$ may also be of interest.

\medskip The sampling averages will, of course, converge when $\Psi$ is mean ergodic.
Unitary evolutions in Hilbert space, for example, will guarantee converging
sampling averages. But sampling averages may possibly also converge on evolutions that are not mean ergodic. The sampling convergence on finitary evolutions is characterized as follows.

\begin{theorem}[Sampling Theorem]\label{t.sampling-theorem} Let $\Psi=(\psi,s)$ be an
arbitrary finitary evolution and $f:{{\mathcal{U}}} \to \cF$ a sampling
function. Then the following statements are equivalent:
\begin{enumerate}
\item[(i)] The sampling averages $\overline{f}_{t}$ converge.
\item[(ii)] The sampling values $f_{t}$ are bounded in norm.
\end{enumerate}
\end{theorem}
Again, we defer the \emph{proof} of Theorem~\ref{t.sampling-theorem} to the Appendix. Choosing $\cF = \cU$ and $f=I$ as the identity operator, we immediately note:

\begin{corollary}\label{c.finitary-mean-ergodic}
 Let $\Psi=(\psi,s)$ be an arbitrary finitary evolution in $\cU$. Then
$$
\mbox{$\Psi$ is mean ergodic} \quad\Longleftrightarrow\quad \mbox{$\Psi$ is stable.}
$$
\qed
\end{corollary}

\section{Observables and Markov chains}\label{sec:observables}

\label{sec:M-chains} Let $A$ be a finite or countable set. An
\emph{observable} with range $A$ on the evolution $\Psi=(\psi,s)$ in
${{\mathcal{U}}}$ is a collection $X=\{\chi_{a}\mid  a\in A\}$ of continuous linear functionals $\chi_{a}$
such that the $\chi_{a}(s^{(t)})$ are real numbers with the property
\[
p_{a}^{(t)}=\chi_{a}(s^{(t)})\geq0\;\;\forall a\in A\quad\mbox {and}\quad
\sum_{a\in A}p_{a}^{(t)}=1.
\]
We think of $X$ as producing the event $\{X_{t}=a\}$ at time $t$ with
probability
\[
\mbox {\rm Pr}\{X_{t}=a\}=p_{a}^{(t)}.
\]
So the observable $X$ yields a sequence $(X_{t})$ of stochastic variables
$X_{t}$ with probability distributions $p^{(t)}$. We call $(X_{t})$ a
\emph{(generalized) Markov chain} on $A$.

\medskip
\begin{Bemerkung} The probability distributions $p^{(t)}_a$ of $X$ may be viewed  as
 ''stochastic kernel'' of $X$ and thus generalize the idea of a kernel
of classical Markov chain theory (see, \emph{e.g.} \cite{FellerII,Hernandez-Lassere}). A
Markov chain in our sense, however, does not need to be a stochastic process nor does its
stochastic kernel need to reflect any conditional probabilities with respect to state transitions.
\end{Bemerkung}

\medskip Related to the (generalized) Markov chain $(X_{t})$ is are
\emph{(statistical) sampling processes} $(Y_{t}^{a})$ with respect to any
$a\in A$, where
\[
Y_{t}^{a}=\left\{
\begin{array}
[c]{cl}%
1 & \mbox {if $X_t=a$}\\
0 & \mbox {otherwise}.
\end{array}
\right.
\]
$(Y_{t}^{a})$ is a Markov chain on $\Psi$ in its own right with binary
alphabet $\{0,1\}$ and probability distributions
\[
\mbox {\rm Pr}\{Y_{t}^{a}=1\}= p_{a}^{(t)}\quad\mbox {and}\quad
\mbox {\rm Pr}\{Y_{t}^{a}=0\}= 1 -p_a^{(t)}.
\]

Assuming that the evolution $\Psi$ is finitary, for example, Theorem~\ref{t.sampling-theorem} implies that
the expectation of the observation averages
\[
E\left(  \frac{1}{t}\sum_{m=1}^{t}Y_{m}^{a}\right)  =\frac{1}{t}\sum_{m=1}%
^{t}E(Y_{m}^{a})=\frac{1}{t}\sum_{m=1}^{t}\mbox {\rm Pr}\{Y_{m}^{a}=1\}
\]
converges to a definite limit $\ov{p}^{(\infty)}_a\geq 0$ since the numbers $\mbox {\rm Pr}\{Y_{t}^{a}=1\}
$ are bounded. In view of
$$
 \sum_{a\in A} p_a^{(m)} =1  \quad\mbox{for all $m\geq 1$},
$$
we conclude that $\ov{p}^{(\infty)}= \{\ov{p}^{(\infty)}_a|a\in A\}$ is a probability distribution on $A$. It is in that sense that we refer to $\ov{p}^{(\infty)}$ as the \emph{limit distribution} of the Markov chain $(X_t)$.

\section{Examples}\label{sec:examples}

\subsection{Evolutions of stochastic processes}\label{sec:stochastic-processes}
Let $(X_t)$ be a discrete stochastic process that takes values in the alphabet $A$. Without loss of generality, we assume that $A$ is binary, say $A=\{0,1\}$. As usual, we denote the set of all finite length words over $A$ as
$$
      A^* = \bigcup_{n=0}^\infty A^n,
$$
where $A^0 =\{\Box\}$ and $\Box$ is the empty word. For any $v =v_1v_2\ldots v_n\in A^n$, $|v| = n$ is the length of $v$. Recall that $A^*$ is a semigroup with neutral element $\Box$ under the \emph{concatenation} operation
$$
    (v_1\ldots v_n)(w_1\ldots w_k)= v_1\ldots v_nw_1\ldots w_k.
$$
Moreover, we set
$$
    p(w_1\ldots w_k|v_1\ldots v_n) = \Pr\{X_{n+1}=w_1,\ldots X_{n+k}= w_k | X_1=v_1,\ldots, X_k = v_n\}
$$
and $p(v) = p(v|\Box)$. If $v\in A^t$ has been produced, the process is in a state that is described by the \emph{prediction vector} $P^v$ with the components
$$
  P^v_w = p(w|v) \quad(w\in A^*).
$$
The expected coordinate values of the next prediction vector are then given by the components of the vector
$$
      \psi P^v = p(0|v)P^{v0} + p(1|v)P^{v1}.
$$
Binomial expansion, therefore, immediately shows that the expected prediction vector at time $t$ is given by
$$
    \sum_{v\in A^t} p(v)P^v = \psi^t P^\Box.
$$
Since the components of the prediction vectors $P^v$ are bounded, they generate a normed vector space $\cP$ with respect to the supremum norm
$$
    \|g\|_\infty = \sup_{w\in A^*} |g_w|.
$$
Moreover, is it not difficult to see that $\psi$ extends to a unique linear operator on $\cP$. The evolution $\Psi=(\psi, P^\Box)$ is said to be the \emph{evolution} of the stochastic process $(X_t)$.

\medskip
For any $a\in A$, one has
$$
(\psi^t P^\Box)_a = \sum_{v\in A^t} p(v)P^v_a = \Pr\{X_{t+1} = a\}.
$$
Since coordinate projections are continuous linear functionals, we find that observations on $(X_t)$ yield observables in the sense of Section~\ref{sec:observables},  which allows us to view the stochastic process $(X_t)$ as a (generalized) Markov chain on $A$.

\medskip\begin{Bemerkung}
The evolution of stochastic processes was first studied by Faigle and Sch\"on\-huth~\cite{FS}, to which be refer for further details. The stochastic evolution model generalizes earlier linear models for the analysis of Markov
type stochastic processes (\emph{e.g.},  Gilbert~\cite{Gilbert}, Dharmadhikari~\cite{Dharma65}, Heller~\cite{Heller} and Jaeger~\cite{Jaeger}).
\end{Bemerkung}

\subsection{Finite-dimensional evolutions}\label{sec:finite-dimensional-evolutions}
Finite-dimensional evolution models are of particular interest in applications. An evolution $\Psi=(\psi, s)$ in $\R^n$ admits an ($n\times n$)-matrix $M$ such that $\psi x = Mx$ holds for all $x\in\R^n$. As observed in \cite{FS}, $\Psi$ is mean ergodic exactly when $\Psi$ is stable (\emph{cf.} Corollary~\ref{c.finitary-mean-ergodic}).

\medskip
Letting $N=\{1,\ldots,n\}$ and assuming that $s$ and all columns of $M$ are probability distributions on $N$, $\Psi$ is clearly stable and hence mean ergodic. The Markov chain relative to $\{\chi_i|i\in N\}$, where the $\chi_i$ are the projections onto the $n$ components of $x\in \R^n$, yields the well-known model of a \emph{random walk} on $N$ with the \emph{transition matrix} $M$.  Considering any map $X:N\to A$ into some alphabet $A$, induced Markov chain with the ''kernel'' functionals
$$
    \chi_a(x) = \sum_{X(i) = a} x_i \quad(x=(x_1,\ldots,x_n)\in \R^n)
$$
is classically known as a \emph{hidden Markov chain} on $A$. Hidden Markov models have proved very useful in practical applications \footnote{see, \emph{e.g.}, \cite{Choi-et-al,Elliot,Vidyasagar}}.

\medskip
It is important to note, however, that even finite-dimensional Markov chains in the general sense of Section~\ref{sec:M-chains} are not necessarily stochastic processes (see Example~\ref{ex.Bell} in Section~\ref{sec:Bell} below). Also quantum random walks (Section~\ref{sec:Q-random-walks} are not necessarily stochastic processes.

\subsection{Quantum statistics}

\label{sec:Q-statistics} Let ${{\mathcal{H}}}$ be a complex Hilbert space of
dimension $|N|$, where $N=\{1,\ldots,n\}$ is finite or $N={{\mathbb{N}}}$,
with inner product $\langle x|y\rangle$. Let ${{\mathcal{B}}}={{\mathcal{B}}%
}({{\mathcal{H}}})$ be the normed complex vector space of all continuous
(linear) operators (\emph{i.e.}, operators $T$ with norm $\|T\|<\infty$). We
single out the set of normalized \emph{wave functions}\footnote{or \emph{qbits}
in the terminology of quantum computing~\cite{Nielsen-Chuang} if $\dim\cH<\infty$}
\[
{{\mathcal{W}}}=\{s\in{{\mathcal{H}}}\mid\|s\|^{2}=\langle s|s\rangle=1\}.
\]
Any $s\in{{\mathcal{W}}}$ gives rise to a (projection) operator $P_{s}%
\in{{\mathcal{B}}}$, where
\[
P_{s}u=\langle u|s\rangle s\quad\mbox {and hence}\quad P_{s}^{2}=P_{s}.
\]
The element $s\in{{\mathcal{W}}}$ furthermore defines a (linear) \emph{trace
functional} $\tau_{s}:{{\mathcal{B}}}\to{{\mathbb{C}}}$ \textit{via}
\[
\tau_{s}T=\langle Ts|s\rangle\quad\mbox {for all $T\in {{\mathcal{B}}}$}.
\]
The trace functional is nonnegative real-valued on every projection operator
$P_{e}$:
\begin{equation}
\label{eq.trace}\tau_{s}P_{e}=\langle P_{e}s|s\rangle=\langle \langle s|e \rangle e|s\rangle=\langle s|e\rangle
\langle e|s\rangle=|\langle s|e\rangle|^{2}\;\in{{\mathbb{R}}}_{+}.
\end{equation}
We therefore obtain from any orthonormal basis $\{e_{i}\mid i\in N\}$ a
probability distribution on $N$ with coefficients $p_{i}=\tau_{s}P_{e_{i}}$:
\begin{equation}
\label{eq.Q-probabilities}p_{i}=| \langle s|e_{i} \rangle |^{2}\geq0\quad\mbox {and}\quad
\sum_{i\in N}p_{i}=\sum_{i\in N}| \langle s|e_{i} \rangle |^{2}=\|s\|^{2}=1.
\end{equation}

\medskip Switching viewpoints, one finds that the collection $X=\{\tau_{e_{i}%
}\mid i\in N\}$ of trace functionals $\tau_{e_{i}}$ yields an observable (in
the sense of Section~\ref{sec:M-chains}) for every evolution $\Psi$ of
projection operators $P_{s}$ in the operator space ${{\mathcal{B}}}$.

\medskip Consider, for example, the Schr\"odinger evolution $\Phi
_{U}=\{s^{(t)}=U^{t}s\mid t\geq0\}$ of the wave function $s\in{{\mathcal{H}}}$
with $\|s\|=1$ relative to the unitary operator $U$. $U$ induces the linear
transformation
\[
T\mapsto UTU^{*}
\]
on ${{\mathcal{B}}}$. Notice that
\begin{equation}
P_sU^{*}(u)=\langle U^{*}(u) | s \rangle s = \langle u | Us \rangle U^{*}U s = U^{*}P_{Us} u
\end{equation}
i.e. $P_{s}U^{*}=U^{*}P_{Us}$ and hence $U%
P_{s}U^{*}=P_{Us}$ holds. So $\Phi_{U}$ has the companion evolution
\[
\Psi_{U}=\{P_{s^{(t)}}=(UP_{s}U^{*})^{t}\mid t\geq0\}
\]
of associated projection operators in ${{\mathcal{B}}}$ that can be observed
under $X$. Slightly more generally, the \emph{states} of a quantum system are
thought to be described by \emph{densities}, \emph{i.e.}, operators $D$ of the
form
\begin{equation}
\label{eq.density}D=\sum_{i\in N}\lambda_{i}P_{e_{i}},
\end{equation}
where $\{e_{i}\mid i\in N\}$ is an orthonormal basis of ${{\mathcal{H}}}$ and
$\{\lambda_{i}\mid i\in N\}$ a (real) probability distribution on $N$.

\subsubsection{Quantum measurements}

\label{sec:Q-measurement} In the standard interpretation of quantum mechanics,
a \emph{measurement} is represented by an operator $M\in{{\mathcal{B}}}$ of
the form
\begin{equation}
\label{eq.Q-measurement}M=\sum_{i\in N}\lambda_{i}P_{e_{i}}%
\end{equation}
where $\{e_{i}\mid i\in N\}$ is an orthonormal basis, and the $\lambda_{i}$
are real (but not necessarily nonnegative) numbers and are the eigenvalues of
$M$. Let $\Lambda$ be the set of different eigenvalues.

\medskip When a quantum system is in the state $P_{s}$ which is implied by the
wave function $s\in{{\mathcal{W}}}$, the measurement is expected to produce
the numerical value
\begin{equation}
E_{M}(s)=\sum_{i\in N}\lambda_{i}\tau_{s}P_{e_{i}}=\sum_{i\in N}\lambda
_{i}p_{i},\label{Q-expectation}%
\end{equation}
where the $p_{i}$ are the probabilities as in (\ref{eq.Q-probabilities}). So
the measurement comes down to the application of the $\Lambda$-valued
observation variable $X$ that takes on a particular value $\lambda$ with
probability
\[
\mbox {\rm Pr}\{X=\lambda\}=\sum_{i:\lambda_{i}=\lambda}\tau_{s}P_{e_{i}}%
\]
and has the expectation
\[
E_{X}(s)=\int_{{\mathbb{R}}}xdp=\sum_{\lambda\in\Lambda}\lambda\mbox
{\rm Pr}\{X=\lambda\}=E_{M}(s).
\]

\medskip In the finite-dimensional case $\operatorname{dim}{{\mathcal{H}}%
}=n<\infty$, the operators $M$ of the form (\ref{eq.Q-measurement}) are
precisely the self-adjoint operators and
\[
E_{M}(s)=\mbox {\rm tr}(MP_{s})
\]
is the usual trace of the product operator $MP_{s}$. If one restricts
attention to Schr\"odinger evolutions,
our quantum statistical model above becomes the
\emph{quantum statistical inference} model proposed by
Barndorff-Nielsen, Gill and Jupp~\cite{B-NGJ}.

\medskip From a mathematical point of view, of course, there is no reason to
restrict statistical inference theory to the ana\-lysis of Schr\"odinger
evolutions. In the same way classical Markov chains generalize to hidden
Markov chains, observable operator models \textit{etc.}, or more general
evolutions in ${{\mathcal{H}}}$ or ${{\mathcal{B}}}({{\mathcal{H}}})$ may be
of interest as well. The statistics of such evolutions can be analyzed in the
same way.

\subsection{Quantum random walks}\label{sec:Q-random-walks}

"Quantum random walks" and "quantum Markov chains" as generalizations of the
classical models to the quantum model have received considerable recent
interest. The models proposed in the literature are typically derived from
Schr\"odinger type evolutions relative to a set $N$. The resulting random walk
is then a particular $N$-valued Markov process in the sense of
Section~\ref{sec:M-chains}.

\medskip As an illustration, we outline a generalization of
Gudder's~\cite{Gudder08} model relative to the set $N=\{1,\ldots,n\}$. Let
$d\geq1$ be some integer parameter and ${{\mathbb{H}}}_{d}$ the real vector
space of all self-adjoint $d\times d$ matrices $C$ with coefficients
$C_{ij}\in{{\mathbb{C}}}$. Define a \emph{state} to be a collection
$S=\{S_{i}\in{{\mathbb{H}}}_{d}\mid i\in N\}$ of self-adjoint matrices $S_{i}$
with nonnegative eigenvalues such that
\[
\mbox {\rm tr}(S)=\sum_{i\in N}\mbox {\rm tr}(S_{i})=1.
\]
Assume to be further given a set $E=\{\epsilon_{ij}\mid i,j\in N\}$ of
operators on ${{\mathbb{H}}}_{d}$ that map densities onto densities. Consider
a process that starts from a state $S$ and iteratively effects state
transitions as follows:
\[
S^{(t)}\mapsto S^{(t+1)}\quad\mbox {with}\quad S_{i}^{(t+1)}=\sum_{j\in
N}\epsilon_{ij}(S_{j}^{(t)})\quad(i\in N).
\]
This process induces an evolution $(S^{(t)})$ in the universe ${{\mathcal{U}}%
}={{\mathbb{H}}}_{d}^{n}$ with evolution matrix $M$, say. Let $\Pi_{i}$ be the
projector that maps $T\in{{\mathbb{H}}}_{d}^{n}$ onto its $i$th coordinate
$T_{i}\in{{\mathbb{H}}}_{d}$ and consider the set of linear operators
\[
{{\mathcal{M}}}=\{M^{(i)}=\Pi_{i}M\mid i\in N\}.
\]
${{\mathcal{M}}}$ induces an $N$-valued stochastic process $(X_{t})$ with
distribution
\[
\mbox {\rm Pr}\{X_{1}=i_{1},\ldots,X_{t}=i_{t}\}=\mbox {\rm tr}(M^{i_{t}%
}(M^{i_{t-1}}(\ldots(M^{(i_{1})}S)\ldots)),
\]
which constitutes a quantum analog of a classical random walk on $N$.

\medskip
\begin{Bemerkung} The quantum random walk model proposed by Aharonov
\textit{et al.}~\cite{Aharonov-et-al} (see also \cite{Faigle11,Kempe,Portugal-et-al-15,Szegedy04})
follows from the present approach by specializing the quantum evolution and observation further. For a
generalization of the classical Metropolis random walk into a quantum context,
see, \emph{e.g.}, Temme \textit{et al.}~\cite{Temme11}.
\end{Bemerkung}

\section{Joint observations}

\label{sec:joint-observations} We say that the  $k$ observables $X^{(1)},\ldots,X^{(k)}$
with alphabets $A_{1},\ldots,A_{k}$ are \emph{jointly observable}
on the evolution $\Psi=(\psi,s)$ if there exists an observable $X$ for $\Psi$ with
alphabet $A=A_{1}\times\ldots\times A_{k}$ such for all $j=1,\ldots,k$,
$a_{j}\in A_{j}$ and $t\geq0$,
\begin{equation}
\label{eq.compatibility}\mbox {\rm Pr}\{X_{t}^{(j)}=a_{j}\}=\displaystyle \sum
_{(a_{1},\ldots,a_{j},\ldots,a_{k})\in A}\mbox {\rm Pr}\{X_{t}=(a_{1}%
,\ldots,a_{j},\ldots,a_{k})\},
\end{equation}
which means that $X^{(j)}$ is the $j$th marginal of $X$.

\medskip It is clear that joint observability of $\{X^{(1)},\ldots,X^{(k)}\} $
implies joint observability for any subset $\{X^{i_{1}},\ldots,X^{i_{r}}\}$.
In particular, sums and products of jointly observable variables are
observable and expected values, covariances \textit{etc.} are well-defined.

\subsection{Heisenberg uncertainty}

\label{sec:Heisenberg} We illustrate the concept of jointly observables with
the example of an important measurement issue in the standard model of quantum
theory. Let $X$ and $Y$ be two observation variables associated with two
quantum measurements as in Section~\ref{sec:Q-measurement}. So there are
representative self-adjoint operators
\begin{equation}
\label{eq.operator-representations}A=\sum_{i\in N}\lambda_{i}P_{e_{i}}%
\quad\mbox {and}\quad B=\sum_{j\in N}\lambda_{j}^{\prime}P_{f_{j}}%
\end{equation}
relative to orthonormal bases $\{e_{i}\mid i\in N\}$ and $\{f_{j}\mid j\in
N\}$ of a complex Hilbert space ${{\mathcal{H}}}$. Let $\Lambda$ and
$\Lambda^{\prime}$ be the ranges of $X$ and $Y$.

\medskip If $X$ and $Y$ are jointly observable in this measurement model
relative to the wave function $s$, there is an observable $Z$ with marginals
$X$ and $Y$ and probability distribution
\[
\mbox {\rm Pr}\{Z=(\lambda,\lambda^{\prime})\}=\sum_{k \in \zeta ^{ -1}\left (\lambda  ,\lambda ^{ \prime }\right )}\tau _{s}P_{g_{k}}
\]

for a suitable map $\zeta:N\to\Lambda\times\Lambda^{\prime}$ and orthonormal
basis $\{g_{k}\mid k\in N\}$. $Z$ admits the operator representation
\[
C=\sum_{k\in N}\mu_{k}P_{g_{k}}\quad
\mbox {with\; $\mu _k=\lambda \cdot \lambda '$
if $\zeta (k)=(\lambda ,\lambda ')$}.
\]
Moreover, we have operator representations for $X$ and $Y$ with respect to the
common basis $\{g_{k}\}$:
\[
\tilde{A}=\sum_{k\in N}\lambda_{k}P_{g_{k}}\quad\mbox {and}\quad\tilde{B}%
=\sum_{k\in N}\lambda_{k}^{\prime}P_{g_{k}}\;\;\mbox {where}\;(\lambda
_{k},\lambda_{k}^{\prime})=\zeta_{k}.
\]
Hence $X$ and $Y$ are seen to satisfy the Heisenberg commutativity condition
for observational compatibility:
\begin{equation}
\label{eq.Heisenberg-condition}\tilde{A}\tilde{B}=C=\tilde{B}\tilde{A}.
\end{equation}
Conversely, if $X$ and $Y$ admit operator representations with respect to a
common basis, it is clear that $X$ and $Y$ are jointly observable relative to
every Schr\"o\-d\-inger evolution.

\medskip
\begin{theorem}\label{t.quantum-observability} Quantum measurements $X$ and $Y$ are
jointly observable on Schr\"o\-dinger evolutions if and only if $X$ and $Y$
admit operator representations with respect to a common orthonormal basis.
\hfill$\blacksquare$
\end{theorem}

\medskip Say that a quantum measurement $X$ has the \emph{Riesz} property  if it admits
a finite or countable orthonormal basis of eigenvectors $e_i$ with eigenvalues $\lambda_i$ such that each
eigenspace $E_\lambda$ is finite-dimensional and
\begin{equation}\label{eq.q-Riesz}
     X = \sum_{i}\lambda_i P_{e_i}.
\end{equation}

\medskip
\begin{corollary}\label{c.Riesz-Heisenberg} Riesz quantum measurements $X$ and $Y $
are jointly observable on Schr\"o\-dinger evolutions if and only if they admit
operator representations by commuting operators.
\end{corollary}

\textit{Proof.$\;$}Let $A$ and $B$ be representations of $X$
and $Y$ as in (\ref{eq.q-Riesz}). It remains to show that
$AB=BA$ implies the existence of a common representative orthonormal basis for
$X$ and $Y$.

\medskip Consider an arbitrary eigenvalue $\lambda\in\Lambda$ of $A$ with
eigenspace $E_{\lambda}$. $A$ is $E_{\lambda}$-invariant. Moreover, for any
$v\in E_{\lambda}$ one has
\[
BAv=\lambda Bv=ABv\quad\mbox {and hence}\quad Bv\in E_{\lambda},
\]
\emph{i.e.}, $E_{\lambda}$ is also $B$-invariant. $E_{\lambda}$ is finite-dimensional and
therefore admits an orthonormal basis $G_{\lambda}$ of eigenelements $g_{\lambda}$ of $B$ that are
also eigenelements of $A$. So
\[
G=\bigcup_{\lambda\in\Lambda}G_{\lambda}%
\]
is an orthonormal basis with the desired property. \hfill$\blacksquare$

\medskip In the same way, one finds:

\medskip
\begin{corollary}\label{c.finite-Heisenberg} The Riesz quantum measurements
$X^{1},\ldots,X^{k}$ on Schr\"odinger evolutions in a Hilbert space are
jointly observable if and only if they are pairwise observable. \hfill
$\blacksquare$
\end{corollary}

\subsection{A Bell-type inequality}\label{sec:Bell}

We have seen that Riesz quantum measurements on
Schr\"odinger evolutions are jointly observable if and only if they are
\emph{pairwise} observable (Corollary~\ref{c.finite-Heisenberg}). This
convenient criterion no longer applies to observables on arbitrary Markov
chains (see Example~\ref{ex.Bell} below). We now establish a necessary condition for
joint observability of three observables in the spirit of Bell's~\cite{Bell66}
inequalities in the standard quantum model.

\medskip
\begin{lemma}[Bell inequality]\label{l.Bell} Let $X,Y,Z$ be pairwise observable on
the (arbitrary) evolution $\Psi=(\psi,s)$, each taking values in $\{-1,+1\}$.
Then the inequality
\begin{equation}
\label{eq.Bell}|E_{t}(XY)-E_{t}(YZ)|\;\leq\;1-E_{t}(XZ)\quad\mbox {holds for
all $t\geq 0$,}
\end{equation}
where $E_{t}(XY)$ is the expected value of the product variable $X_{t}Y_{t}$
at time $t$.
\end{lemma}

\textit{Proof.$\;$}Any choice of $x,y,z\in\{-1,+1\}$ satisfies the inequality
\[
|xy-yz|\;\leq\;1-xz.
\]
The probabilities $p_{t}(x,y,z)=\mbox {\rm Pr}\{X_{t}=x,Y_{t}=y,Z_{t}=z\}$ are
nonnegative real numbers that sum up to $1$. So we conclude
\begin{eqnarray*}
|E_{t}(XY)-E_{t}(YZ)| & =&\big |\sum_{x,y,z}(xy-yz)p_{t}(x,y,z)\big |\;\leq
\;\sum_{x,y,z}|xy-yz|p_{t}(x,y,z)\\
& &\leq\sum_{x,y,z}(1-xz)p_{t}(x,y,z)\;=\;1-E_{t}(XZ)\;.
\end{eqnarray*}
\hfill$\blacksquare$

\medskip The inequality (\ref{eq.Bell}) may be violated by observables $X,Y,Z$
that are \emph{pairwise} but not \emph{jointly} observable.

\medskip\begin{ex}\label{ex.Bell} Consider the (stationary) evolution
$\Psi=(\psi,D)$ with $D^{(t)}=D$ in the space ${{\mathbb{H}}}_{5}$ of all
$5\times5$ self-adjoint matrices, where
\[
D=\mbox {\rm diag}(-1/3,1/3,1/3,1/3,1/3)\in{{\mathbb{H}}}_{5}%
\]
and $\mbox {\rm diag}(v)$ denotes the diagonal matrix with diagonal vector $v
$. Let $X,Y,Z$ be the measurements that are induced be the self-adjoint
matrices
\begin{eqnarray*}
A_{X} & =&\mbox {\rm diag}(-1,+1,-1,-1,-1)\\
A_{Y} & =&\mbox {\rm diag}(+1,+1,-1,+1,-1)\\
A_{Z} & =&\mbox {\rm diag}(+1,+1,+1,-1,-1).
\end{eqnarray*}
Notice that $A_{X},A_{Y},A_{Z}$ commute pairwise and thus satisfy the
Heisenberg condition (\ref{eq.Heisenberg-condition}). Moreover, $X,Y,Z$ are
pairwise observable on $\Psi$. The pairs of products have the expectations
\begin{equation}
E(XY)=+1,\;E(YZ)=-1/3,\;E(XZ)=+1\
\end{equation}
and violate the Bell inequality (\ref{eq.Bell}), which shows that $X,Y,Z$ are
not jointly observable on $\Psi$.
\end{ex}

\medskip
\begin{Bemerkung} The experimental results of Aspect \textit{et
al.}\cite{Aspect} suggest that measurements on real world quantum systems may
violate Bell's inequalities. In our Markov setting, these results can be
explained as follows: the experiments were either not carried out with
pairwise commuting observables (and thus subject to Heisenberg uncertainty)
and/or the description of quantum states by ''densities'' with only
\emph{nonnegative} eigenvalues is too restrictive for real world models.
\end{Bemerkung}

\section{Conclusion}
A model for the Markovian statistical analysis of observations on evolving systems has been proposed that separates the evolution of the system states and the observation of system events clearly. This model not only generalizes classical views on homogeneous Markov chains as random sources properly but allows a Markov type analysis of more general observation processes which, in particular, include observations arising from general underlying stochastic processes. This separation of the notions of system evolution and system observation allows us to develop a general theory of joint observability, which has no classical counterpart. It is compatible with the notion of Heisenberg uncertainty relative to Schr\"odinger evolutions but it is not implied by it.

\medskip
Several intriguing questions immediately raise themselves. For example, we do not think that the model of Riesz evolutions is  the most general in which Markovian convergence can be proved. Do our convergence results hold relative to evolution operators $T$ with $\lambda =1$ being just an isolated eigenvalue of $T$? What is the general convergence behavior of quantum densities? Is it true, for instance, that a normal operator $\psi$ of norm $\|\psi\|\leq 1$ on a Hilbert space $\cH$  not only yields a mean ergodic evolution in $\cH$ itself (Theorem~\ref{t.vNeumann-normal-mean}) but also a mean-ergodic evolution of the associated densities?

\section{Appendix: Proofs}

For fundamental notions and facts on linear operators we refer to standard
texts\footnote{\emph{e.g.}, \cite{Conway90,Dowson78}} for
fundamentals on linear operators.

\subsection{Proof of Theorem~\ref{t.vNeumann-normal-mean}}

Recall the spectral representation for a continuous normal operator in its
multiplication form:

\begin{theorem}\label{t.normal-spectral-representation} Let $\psi$ be a continuous operator
on a complex Hilbert space ${{\mathcal{H}}}$ such that $\psi^{\ast}\psi
=\psi\psi^{\ast}$. Then there exists a measure space $(\Omega,\Sigma,\mu)$, an
essentially bounded measurable function $g:\Omega\rightarrow{{\mathbb{C}}}$
and a unitary operator $U:{{\mathcal{H}}}\rightarrow{{\mathcal{L}}}_{\mu}%
^{2}(\Omega)$ such that
\[
\psi=U^{\ast}M_{g}U,
\]
where $M$ is the multiplication operator $M_{g}f=f\cdot g$. Moreover,
\[
\Vert\psi\Vert=\Vert M_{g}\Vert=\Vert g\Vert_{\infty}.
\]
\end{theorem}

\medskip By Theorem~\ref{t.normal-spectral-representation}, we can assume w.l.o.g.:

\begin{itemize}
\item[$\bullet$] ${{\mathcal{H}}}={{\mathcal{L}}}_{\mu}^{2}(\Omega)$ and
$\psi$ is given as multiplication by a bounded measurable function $g$.
\end{itemize}

The stability of the evolution $\Psi=(\psi,s)$ implies that there is a constant $M$ so that $|g^n(\omega)f(\omega)| \leq M$ holds almost everywhere for all positive integers $n$. It follows that a.e.,   $|g(\omega)| \leq 1$ or $f(\omega) = 0$. Hence there is a measurable function $g_1$ so that a.e.,  $|g_1(\omega)|\leq1$ and $g^n(\omega)f(\omega) = g_1^n(\omega)f(\omega)$ . Setting
\[
\overline{\psi}_{t}(f)=\left(  \frac{1}{t}\sum_{m=1}^{t}g_1^{m}\right)  f,
\]
we therefore conclude
\[
|\overline{\psi}_{t}(f)(\omega)|^{2}\leq\left(  \frac{1}{t}\sum_{m=1}%
^{t}|g_1(\omega)|^{m}\right)  ^{2}|f(\omega)|^{2}\leq|f(\omega)|^{2}.
\]
The sequence $(\overline{\psi}_{t}(f)(\omega))_{t \geq 0}$ converges to
\[
\pi(f)\left(  \omega\right)  =\left\{
\begin{array}
[c]{cl}%
f(\omega) & \mbox {if $g_1(\omega )=1$}\\
0 & \mbox {otherwise.}
\end{array}
\right.
\]
If $g_1(\omega)=1$, then $\overline{\psi}_{t}(f)(\omega)=f(\omega)=\pi(f)\left(
\omega\right)  $, and so
\begin{eqnarray*}
\lim_{t\rightarrow\infty}\Vert\overline{\psi}_{t}(f)(\omega)\Vert_{2}  &
=&\lim_{t\rightarrow\infty}\left(  \int_{\Omega}|\overline{\psi}_{t}%
(f)-\pi(f)(\omega)|^{2}d\omega\right)  ^{1/2}\\
& =&\lim_{t\rightarrow\infty}\left(  \int_{g(\omega)\neq1}|\overline{\psi}%
_{t}(f)\left(  \omega\right)  |^{2}d\omega\right)  ^{1/2}.
\end{eqnarray*}
On the set $\{\omega\in\Omega\mid g(\omega)\neq1\}$, the functions
$|\overline{\psi}_{t}(f)(\omega)|^{2}$ converge pointwise to $0$ and are
bounded by the integrable function $|f(\omega)|^{2}$. The theorem of dominated
convergence thus yields
\[
\lim_{t\rightarrow\infty}\left(  \int_{g(\omega)\neq1}|\overline{\psi}%
_{t}(f)\left(  \omega\right)  |^{2}d\omega\right)  ^{1/2}=\left(
\int_{g(\omega)\neq1}\lim_{t\rightarrow\infty}|\overline{\psi}_{t}(f)\left(
\omega\right)  |^{2}\right)  ^{1/2}d\omega=0.
\]
Clearly, $\pi(f)$ is the orthogonal projection of $f$ onto the eigenspace of $\lambda=1$. So stability is
sufficient for mean ergodicity.

\medskip
To see that stability is necessary for mean ergodicity, assume that
$$
\lim_{n\to\infty}\frac{1}{n}\sum _{k =0}^{n -1}T^{k}\left (f\right )\;\;\mbox{exists for the operator}\quad T\left (f\right ) =\int gfd{\mu}.
$$
We would like to show that $\left \vert g\left (x\right )\right \vert  \leq 1$ holds a.e. on the set $\left \{x :f\left (x\right ) \neq 0\right \}$.

\medskip
Assume that the set
\begin{equation}M =\left \{x :f\left (x\right ) \neq 0\text{~}\ \text{and~}\ \text{}\text{}\left \vert g\left (x\right )\right \vert  >1\right \}
\end{equation}has positive measure. Then for some integer $r >1$ the set\begin{equation}M_{r} =\left \{x :f\left (x\right ) \neq 0\text{~}\ \text{and~}r >\ \text{}\text{}\left \vert g\left (x\right )\right \vert  >1 +\frac{1}{r}\right \}
\end{equation}has positive measure and for any $x \in M_{r}$,  we have
\begin{equation}\frac{1}{n}\sum _{k =0}^{n -1}T^{k}\left (f\right ) =\frac{f\left (x\right )}{n}\sum _{k =0}^{n -1}g(x)^{k} =\frac{f\left (x\right )}{n}\frac{1 -g\left (x\right )^{n}}{1 -g\left (x\right )}.
\end{equation}
Observing
\begin{eqnarray*}
&&\left \vert 1 -g\left (x\right )^{n}\right \vert  \geq \left \vert g\left (x\right )\right \vert ^{n} -1 >\left (1 +\frac{1}{r}\right )^{n}-1\\
{\rm and}&&\left \vert 1 -g\left (x\right )\right \vert  \leq 1 +\left \vert g\left (x\right )\right \vert  \leq 1 +r,
\end{eqnarray*}
we thus conclude
\begin{equation}\left \vert \frac{1}{n}\sum _{k =0}^{n -1}T^{k}\left (f\right )\left (x\right )\right \vert  \geq \frac{\left \vert f\left (x\right )\right \vert }{n}\frac{\left (1 +\frac{1}{r}\right )^{n}-1}{1 +r}
\end{equation}
and hence
\begin{eqnarray*}\left \Vert \frac{1}{n}\sum _{k =0}^{n -1}T^{k}\left (f\right )\right \Vert _{2} &=&\frac{1}{n}\left (\int _{X}\left \vert \sum _{k =0}^{n -1}T^{k}\left (f\right )\left (x\right )\right \vert ^{2}\right )^{1/2} \\
 &\geq& \frac{1}{n}\left (\int _{M_{r}}\left \vert \sum _{k =0}^{n -1}T^{k}\left (f\right )\left (x\right )\right \vert ^{2}\right )^{1/2} \\
 &\geq& \frac{\left (1 +\frac{1}{r}\right )^{n}-1}{n\left (1 +r\right )}\left (\int _{M_{r}}\left \vert f\left (x\right )\right \vert ^{2}d\mu \right )^{1/2S}.
\end{eqnarray*}
Since $\int _{M_{r}}\left \vert f\left (x\right )\right \vert d\mu  >0$ and $\lim _{n \rightarrow \infty }\frac{\left (1 +\frac{1}{r}\right )^{n}-1}{n\left (1 +r\right )} =\infty  ,$ it follows that $\left \Vert \frac{1}{n}\sum _{k =0}^{n -1}T^{k}\left (f\right )\right \Vert _{2}$ is unbounded. So the series $\left (\frac{1}{n}\sum _{k =0}^{n -1}T^{k}\left (f\right )\right )_{n >0}$ cannot converge. Consequently,

\medskip
$\left \vert g\left (x\right )^{n}f\left (x\right )\right \vert  \leq \left \vert f\left (x\right )\right \vert $ and hence~$\left \Vert T^{n}\left (f\right )\right \Vert _{2} \leq \left \Vert f\right \Vert _{2}$ holds for all~$n$.

\qed

\subsection{Proof of Theorems~\ref{t.Riesz} and \ref{t.sampling-theorem}}

Throughout this section, let $U$ be a fixed Banach space with a fixed element $s \in U$. We further fix an operator $T : U \rightarrow U$  that is bounded on
$U_s = \lin \{T^n s : n\geq 0 \}$
and denote by $\hat{T}$ its (bounded) extension to $\hat{U}_s$. Without loss of generality, we can therefore assume $\hat{U}_s = U$ and $\hat{T} = T$. $\sigma(T)$ denotes the spectrum of $T$. The \emph{spectral radius} of $T$ is
\[
r(T)=\max\{|\lambda
|\mid\lambda\in\sigma(T)\} = \lim_{t\rightarrow\infty}\Vert T^{t}\Vert^{1/t}\leq \|T\|.
\]

\begin{lemma}\label{l.Riesz-spectrum} If $T$ is Riesz, then $\sigma_1(T) =\{\lambda \in \sigma(T)\mid |\lambda| \geq 1\}$
is a finite set.
\end{lemma}

\Pf $\sigma_{\varepsilon,\delta}(T) = \{\lambda\in \sigma(T) \mid \varepsilon \leq |\lambda|\leq \delta\}$
is a bounded subset of $\C$ for any $\varepsilon,\delta\geq 0$. If $T$ is Riesz, $0$ is the only possible accumulation point of $\sigma(T)$. So $\sigma_{\varepsilon,\delta}(T)$ must be finite for every $\varepsilon >0$. If $T$ is bounded, $\sigma_1(T) = \sigma_{\varepsilon,\delta}(T)$ holds for $\varepsilon = 1$ and $\delta = \|T\|$ and the claim of the Lemma follows.

\qed

\medskip
For any eigenvalue $\lambda$ of $T$ with finite algebraic multiplicity $n_\lambda$, the \emph{Riesz de\-com\-position} of $T$ with respect to $\lambda$ guarantees\footnote{see, \emph{e.g.}, \cite{Dowson78}}:
\begin{itemize}
\item[(R)] $U$ admits the direct sum decomposition $U = N_\lambda \oplus R_\lambda$, where
\begin{enumerate}
\item $N_{\lambda}=\{x\in U\mid(T-\lambda)^{n_{\lambda}}x=0\}$ is $T$-invariant with $\dim N_\lambda <\infty$;
\item $R_{\lambda}=(T-\lambda)^{n_{\lambda}}U$ is $T$-invariant.
\end{enumerate}
\end{itemize}

If $\sigma_1(T)$ is finite set of eigenvalues, repeated application of the Riesz decomposition (R) to some $\lambda\in \sigma_1(T)$ and then to $T:R_{\lambda}\to R_{\lambda}$
\textit{etc.} and the other eigenvalues in $\sigma_{1}(T)$ yields

\begin{lemma}[Riesz decomposition]\label{l.Riesz-decomposition} If $\sigma_1(T)$ is a finite set of eigenvalues of $T$ with finite algebraic multiplicities, then $U$ admits a direct sum decomposition $U=N\oplus W$ into $T$-invariant subspaces $N$ and $W$, where
\[
N=\bigoplus_{\lambda\in\sigma_{1}(T)}N_{\lambda} \quad\mbox{and}\quad \dim N <\infty.
\]
Moreover, $|\lambda| < 1$ holds for all eigenvalues $\lambda$ of the restriction of $T$ to $W$.
\qed
\end{lemma}

\medskip
The decomposition (R) implies that Riesz evolutions are finitary.

\begin{proposition}\label{p.1} Assume that $T$ is a Riesz operator with decomposition $U = N\oplus W$ into $T$-invariant subspaces $N$ and $W$ such that $\dim N<\infty$ and the restriction of $T$ to $W$ has no eigenvalue in $\sigma_1(T)$. Then the Riesz evolution $(T,s)$ is equivalent to the finite-dimensional evolution $(T, s_N$), where $s_N\in N$ is such that  $s=s_N +s_W$ holds for some $s_W\in W$.
\end{proposition}

\Pf In view of $T^m s = T^m s_N + T^m s_W$ for all $m\geq 0$, it suffices to establish the claim
$$
   \lim_{n\to\infty} T^n s_W = 0.
$$
Since $\sigma_{\varepsilon, 1}(T)$ is a finite set for any $\varepsilon >0$, the spectral radius $r_W(T)$ of $T$ on $W$
must satisfy  $r_W(T)< 1$. For clarity of notation, let $T_W$ be the restriction of $T$ to $W$ and
choose $n_{0}$ so large that $\Vert T_W^{n}\Vert^{1/n}\leq
r<1$ holds for all $n\geq n_{0}$. Then one has $\Vert T_W^{n}\Vert\leq r^{n}$
and thus concludes
$$
\lim_{n\rightarrow\infty}\left\Vert T_W^{n}s_W\right\Vert =0.
$$
\hfill$\blacksquare$

\medskip
The proof of Theorem~\ref{t.Riesz} is now immediate: The Riesz evolution $(\psi,s)$ is equivalent to the finite-dimensional evolution $(\psi,s_N)$. The ergodic properties stated in Theorem~\ref{t.Riesz} are directly obtained by applying Proposition~\ref{p.FS} to $(\psi,s_N)$.

\qed

\medskip
For the proof of the sampling theorem (Theorem~\ref{t.sampling-theorem}), let $(Q, x)$ be a finite-dimensional evolution that is equivalent to $(T,s)$. So we have $\|T^ns-Q^nx\|\to 0$ and hence
$$
\lim_{n\to\infty} \|f(T^ns) -f(Q^n x)\|  = \lim_{n\to\infty} \|f(T^ns -Q^nx)\| = 0
$$
since the sampling function $f$ is continuous. This implies
$$
\lim_{t\to\infty}\frac1t\sum_{m=1}^t \|f(T^ms) -f(Q^m x)\| =0,
$$
\emph{i.e.}, the $f$-sample averages converge on $(T,s)$ exactly when they converge on $(Q,x)$. It is furthermore clear that $(T,s)$ is stable exactly when $(Q,x)$ is stable. For the proof, we can therefore assume without loss of generality that already $(T,s)$ is finite-dimensional and hence the sample space $\cF = f(U_s)$ is finite-dimensional.

\medskip
Passing to coordinates, we may thus assume: $U_s=\C^n$ and $\cF =\C^k$. Since $f:\C^n\to \C^k$ is bounded on $(T,s)$ if and only if each component functional $f_j$ of $f$ is bounded on $(T,s)$, it suffices to consider the $1$-dimensional case $k=1$.

\medskip
With respect to the chosen coordinatization, $T$ is an $n\times n$ matrix, $s$ a column vector and $f$ a row vector of dimension $n$. Assume first that the sequence $(fT^ts)$ is bounded. Then
$$
\mbox{$\;(fT^tu)\;$ is bounded for every $u\in U_s = \lin\{T^ts\mid t=0,1,\ldots\} =\C^n$.}
$$
It follows that the sequence $(fT^t)$ of $n$-dimensional row vectors constitutes a bounded evolution. (The choice of $u$ as the unit vector $e_i$ in $\C$ shows that the $i$th coordinate of the evolution is bounded.) In view of Proposition~\ref{p.FS} (Section~\ref{sec:equivalent-evolutions}), this evolution is mean-ergodic. Consequently, the boundedness of $(fT^ts)$implies the existence of
$$
  \lim_{t\to\infty} \frac1t\sum_{m=1}^t fT^m \quad\mbox{and hence of}\quad   \ov{f}_\infty= \lim_{t\to\infty} \frac1t\sum_{m=1}^t fT^ms.
$$
To prove the converse implication, assume that $\ov{f}_\infty$ exists. Then
$$
\mbox{$\;\D\frac1t\sum_{m=1}^t \lim_{t\to\infty} fT^mu \;$ exists for every $u\in U_s = \lin\{T^ts\mid t=0,1,\ldots\} =\C^n$.}
$$
It follows that the evolution $(fT^t)$ of row vectors is mean ergodic and hence, again by  Proposition~\ref{p.FS}, stable, \emph{i.e.}, there is some constant $c\in\R$ such that
$$
  |fT^ts| \leq   \|fT^t\|\cdot\|s\| \leq c\|s\| < \infty \quad\mbox{for all $t\geq 0$,}
$$
which establishes the claim of the sampling theorem.

\qed

\end{document}